\documentclass[12pt,twoside]{article}

\usepackage{amssymb,amsmath}

\newtheorem{theorem}{Theorem}[section]

\newtheorem{lemma}[theorem]{Lemma}

\numberwithin{equation}{section}

\newcommand{\N}{\mathbb{N}}

\newcommand{\C}{\mathbb{C}}
\newcommand{\dis}{\displaystyle}

\newcommand{\db}{\rule[.05in]{.09in}{.10in}} %darkbox

\begin{document}

\markboth{Igor E. Pritsker}{Convergence of Julia polynomials}%

\title{Convergence of Julia polynomials}%
\author{Igor E. Pritsker\thanks{Research supported in part
by the National Security Agency under Grant No. MDA904-03-1-0081.}}%

\date{}%

\maketitle

\begin{abstract}

We study the approximation of conformal mappings with the
polynomials defined by Keldysh and Lavrentiev from an extremal
problem considered by Julia. These polynomials converge uniformly
on the closure of any Smirnov domain to the conformal mapping of
this domain onto a disk. We prove estimates for the rate of such
convergence on domains with piecewise analytic boundaries,
expressed through the smallest exterior angles at the boundary.
\\
\\
{\bf 2000 MCS.} Primary 30C40, 30E10; Secondary 41A10, 30C30.

\end{abstract}

% ----------------------------------------------------------------

\section{Extremal problems in Smirnov spaces and the associated polynomials}

Let $G$ be a Jordan domain in the complex plane, bounded by a
rectifiable curve $L$ of length $l$. We consider the Smirnov
spaces $E_p(G), \ 0<p<\infty,$ of analytic functions in $G$, whose
boundary values satisfy
\[
\|f\|_p := \left( \int_{L} |f(z)|^p | dz | \right)^{1/p}< \infty
\]
(see Duren \cite{Du}, Smirnov and Lebedev \cite{SL}). Julia
\cite{Ju} studied the extremal problem
\begin{align} \label{1.1}
\inf \{ \|f\|_p : f\in E_p(G),\ f(\zeta)=1) \},
\end{align}
where $\zeta\in G$ is a fixed point. He showed that the above
infimum is attained by the function $(\phi')^{1/p}$, where $\phi$
is the conformal mapping of $G$ onto a disk $D_R:=\{w:|w|<R\}$,
normalized by the conditions $\phi(\zeta)=0,\ \phi'(\zeta)=1.$ It
is known that one can define sequences of polynomials associated
with many extremal problems. Thus Keldysh and Lavrentiev (cf.
\cite{KL35}, \cite{Ke} and \cite{KL37}) considered the polynomials
$Q_{n,p}(z)$ that minimize \eqref{1.1} among all polynomials
$P_n(z)$ of degree $n$, such that $P(\zeta)=1.$ These extremal
polynomials were attributed to Julia by Keldysh and Lavrentiev in
\cite{KL35}. Perhaps, it is more appropriate to call them
Julia-Keldysh-Lavrentiev polynomials. The goal of such
construction is that, provided polynomials are dense in $E_p(G)$,
it would furnish the approximation to $(\phi')^{1/p}$ by $Q_{n,p}$
as $n\to\infty.$ Keldysh and Lavrentiev \cite{KL37} characterized
a class of domains for which these desired properties hold true.
Their results are summarized below.
\\ \\
\noindent{\bf Theorem KL} {\it Let $\psi=\phi^{-1}.$ The following
conditions are equivalent:
\begin{enumerate}
\item $\log|\psi'|$ is representable in $D_R$ by the Poisson
integral of its boundary values

\item $\dis\lim_{n\to\infty} \|(\phi')^{1/p}-Q_{n,p}\|_p = 0$

\item $Q_{n,p}$ converge to $(\phi')^{1/p}$ locally uniformly in
$G$

\item Polynomials are dense in norm in $E_p(G)$
\end{enumerate}
}

The domains with property {\it 4} were later named after Smirnov.
Although there is no complete geometric description of such
domains, many standard classes of domains possess the Smirnov
polynomial density property. The widest currently known class is
probably that of Ahlfors-regular domains (cf. Pommerenke
\cite[Chap. 7]{Po}). They are defined by the condition that there
exists a constant $C>0$ such that
\[
|L\cap D_r| \le C r,
\]
for any disk $D_r$ of radius $r$, where $|\cdot|$ denotes the
arclength measure on $L=\partial G.$

Since $Q_{n,p}$ converge locally uniformly to $(\phi')^{1/p}$  in
Smirnov domains, they have no zeros in any fixed compact subset of
$G$ for large $n$, by Hurwitz's theorem. Hence the following
functions
\begin{align} \label{1.2}
J_{n,p}(z):=\int_{\zeta}^z (Q_{n,p}(t))^p\,dt,\qquad z\in G,
\end{align}
are well defined for $p\in(0,\infty)$, when $n$ is large. Of
course, they are well defined polynomials for any $n\in\N,$ if
$p\in\N.$ Let $\|\cdot\|_{\infty}$ denote the uniform (sup) norm
on $\overline G.$
\begin{theorem} \label{thm1.1}
If $p\in\N$ then
\[
\|\phi-J_{n,p}\|_{\infty} \le \frac{1}{2}\left((2\pi R)^{1/p} +
l^{1/p} \right)^{p-1} \left\|(\phi')^{1/p}-Q_{n,p}\right\|_p.
\]
Hence $\|\phi-J_{n,p}\|_{\infty} \to 0\mbox{ as } n\to\infty$,
when $G$ is a Smirnov domain.
\end{theorem}
For $p=1$, the uniform convergence of $J_{n,1}$ on $\overline G$
was already observed by Keldysh and Lavrentiev in \cite{KL37}, but
without any estimate. The case $p=2$ is of special interest
because of its close connections with the Szeg\H{o} kernel and
orthogonal polynomials. Ahlfors \cite{Ah} considered these
polynomials for the numerical approximation of conformal mappings,
and developed interesting representations for $Q_{n,2}$ via
iterated integrals of Vandermonde determinants.  For $p=2$, a
similar result to Theorem \ref{thm1.1} was proved by Warschawski
\cite{Wa}. Gaier \cite[pp. 130-131]{Ga64} gave the estimates of
the uniform convergence rates for some smooth domains, based on
the results of Rosenbloom and Warschawski \cite{RW}. The explicit
rates of convergence for domains with corners are stated in the
following section. Note that the case $p=2$ was already considered
in \cite{Pr}, but with a somewhat different normalization for the
polynomials.

A natural analogue of the polynomials $Q_{n,p}$ is given by the
best approximating polynomials $\tilde Q_{n,p}$ to $(\phi')^{1/p}$
in $E_p(G)$:
\begin{align} \label{1.3}
\left\|(\phi')^{1/p} - \tilde Q_{n,p}\right\|_p =\inf \left\{
\left\|(\phi')^{1/p}-P_n\right\|_p : P_n(\zeta)=1 \right\},
\end{align}
where the inf is taken over the polynomials $P_n$ of degree $n$.
Following the same convention as for \eqref{1.2}, we define the
functions
\begin{align} \label{1.4}
\tilde J_{n,p}(z):=\int_{\zeta}^z (\tilde Q_{n,p}(t))^p\,dt.
\end{align}
It can be readily seen that $\tilde Q_{n,2} \equiv  Q_{n,2}$ and
$\tilde J_{n,2} \equiv J_{n,2}$ (for $p=2$), because
\[
\left\|(\phi')^{1/2} - P_n\right\|_2 = \left\|P_n\right\|_2 -
\left\|(\phi')^{1/2}\right\|_2
\]
for any polynomial $P_n$ of degree at most $n,\ P_n(\zeta)=1$, see
\cite[p. 128]{Ga64} and \cite{Ah}. The explicit representation of
$Q_{n,2}$ via the contour orthonormal polynomials $\{ p_n \}_{n
=0}^{\infty}$ in $E_2(G)$ follows from the standard Hilbert space
theory (cf. \cite[Chap. III]{Ga64} and \cite[Chap. 4]{SL}):
\[
Q_{n,2}(z)=\frac{\sum_{k=0}^n \overline{p_k(\zeta)}
p_k(z)}{\sum_{k=0}^n |p_k(\zeta)|^2}, \quad n \in \N.
\]
Thus these polynomials coincide (up to a constant factor) with the
partial sums of the Szeg\H{o} kernel $K(z,\zeta) =
\sum_{k=0}^{\infty} \overline{p_k(\zeta)} p_k(z)$ (cf. Szeg\H{o}
\cite{Sz}). They can be used for the constructive approximation of
the conformal mapping $\phi$, see \cite{Ga64}, \cite{SL} and
\cite{Pr} for the details.

Our interest in $\tilde Q_{n,p}$ and $\tilde J_{n,p}$ is explained
by the fact that one can produce estimates of the convergence
rates for these polynomials. We first state the analogue of
Theorem \ref{thm1.1} in this case.
\begin{theorem} \label{thm1.2}
If $p\in\N$ then
\[
\left\|\phi-\tilde J_{n,p}\right\|_{\infty} \le \frac{1}{2}
\left(3(2\pi R)^{1/p} + l^{1/p} \right)^{p-1}
\left\|(\phi')^{1/p}-\tilde Q_{n,p}\right\|_p.
\]
Hence $\left\|\phi-\tilde J_{n,p}\right\|_{\infty} \to 0\mbox{ as
} n\to\infty$, when $G$ is a Smirnov domain.
\end{theorem}
The rates of convergence quantifying this result are studied in
the next section.

\medskip

We conclude this section by showing that the zeros of the
polynomials $\tilde Q_{n,p}$ are typically dense in $L$. This
indicates that the definition $\tilde J_{n,p}$ in \eqref{1.4}
cannot be extended to $\overline G$ for $p\not\in\N,$ as every
zero generates a branch point. Let $\nu_{n,p}$ and
$\tilde\nu_{n,p}$ be the normalized counting measures for the
zeros of $Q_{n,p}$ and $\tilde Q_{n,p}$ respectively. They are
obtained by placing the point mass $1/n$ at each zero of $Q_{n,p}$
and $\tilde Q_{n,p}$, according to multiplicities. Denote the
equilibrium measure of $\overline G$ (in the sense of logarithmic
potential theory) by $\mu$ \cite{Ra}. The following
Jentzsch-Szeg\H{o} type theorem on the asymptotic zero
distribution is stated in terms of the weak* convergence for
measures.
\begin{theorem} \label{thm1.3}
Let $G$ be a Smirnov domain. Suppose that $\phi$ cannot be
continued as analytic on $\overline G$ function. Then for any
$p\in (0,\infty)$ there exists an infinite subsequence $\tilde
N\subset\N$ such that
\[
\tilde\nu_{n,p} \stackrel{*}{\to} \mu \qquad \mbox{as }
n\to\infty,\ n\in \tilde N.
\]
\end{theorem}
We conjecture that this theorem also holds for the zeros of
$Q_{n,p},$ i.e.,
\[
\nu_{n,p} \stackrel{*}{\to} \mu \qquad \mbox{as } n\to\infty,\
n\in N,
\]
where $N\subset\N$ is an infinite subsequence.  Such asymptotic
behavior is generic for zeros of many extremal polynomials, see
\cite{AB}. It is possible to prove a converse of Theorem
\ref{thm1.3}, i.e., if there is a subsequence of $\tilde\nu_{n,p}$
weakly convergent to $\mu$, then $\phi$ is not analytic on
$\overline G$ (cf. \cite{BSS}).

\section{Convergence in domains with piecewise analytic boundaries}

We consider domains with piecewise analytic boundaries in this
section, which are important in applications. An analytic arc is
defined as the image of a segment under a mapping that is
conformal in an open neighborhood of the segment. Thus a domain
has piecewise analytic boundary if it is bounded by a Jordan curve
consisting of a finite number of analytic arcs. Let $L$ be
piecewise analytic, with the smallest exterior angle $\lambda\pi,\
0<\lambda\le 2,$ at the junction points of the analytic arcs. The
following results contain estimates for the rates of convergence
of $\tilde Q_{n,p}$ in terms of geometric properties of domains.

\begin{theorem} \label{thm2.1}
If $0<\lambda<2$ then
\begin{align} \label{2.1}
\left\|(\phi')^{1/p}-\tilde Q_{n,p}\right\|_p \le C_1 \left\{
\begin{array}{lll}
n^{-\frac{\lambda}{p(2-\lambda)}},\quad &1<p<\infty,\\
n^{-\frac{\lambda}{2-\lambda}} \log n,\quad &p=1,\\
n^{-\frac{\lambda(\lambda-1)}{p(2-\lambda)}-\lambda},\quad
&\frac{1-\lambda}{2-\lambda}<p<1.
\end{array}
\right.
\end{align}
Note that $p\ge 1/2$ works for all $\lambda\in(0,2).$ For
$p\in\N$, we also have
\begin{align} \label{2.2}
\|\phi-\tilde J_{n,p}\|_{\infty} \le C_2 \left\{
\begin{array}{ll}
n^{-\frac{\lambda}{2-\lambda}} \log n,\quad p=1,\\
n^{-\frac{\lambda}{p(2-\lambda)}},\quad p=2,3,\ldots.
\end{array}
\right.
\end{align}
The constants $C_1>0$  and $C_2>0$ are independent of $n\ge 2$.
\end{theorem}

The rates of convergence were previously known only in the case of
$p=2$ for some smooth domains, see \cite[p. 131]{Ga64}. This
theorem is new for any $p\in(0,\infty).$ It is worth noting that
the exponent $\lambda/(2-\lambda)$ for $p=1$ in \eqref{2.1} and
\eqref{2.2} is best possible. Indeed, it is known that this
exponent cannot be improved, in general, for approximation of
$\phi$ in the uniform norm by {\it any} sequence of polynomials
(cf. \cite{Ga92} and \cite{Ga98}). We believe that the exponents
of $n$ are also sharp in \eqref{2.1} and \eqref{2.2} for any
$p\in(1,\infty).$

When all angles at the boundary are the outward pointing cusps, we
can make an even stronger conclusion.

\begin{theorem} \label{thm2.2}
If $\lambda=2$ then for any $p\in (0,\infty)$ there exist $q,r \in
(0,1)$ such that
\begin{align} \label{2.3}
\left\|(\phi')^{1/p}-\tilde Q_{n,p}\right\|_p \le C_3\, q^{n^r},
\qquad n\in\N.
\end{align}
Furthermore, if $p\in\N$ then
\begin{align} \label{2.4}
\|\phi-\tilde J_{n,p}\|_{\infty} \le C_4\, q^{n^r}, \qquad n\in\N.
\end{align}
Here, $C_3>0$  and $C_4>0$ are independent of $n$.
\end{theorem}

It should be mentioned that $r$ cannot be equal to $1$ in the
above theorem, because the geometric rate of convergence implies
that $\phi$ is analytic on $\overline G$ \cite{Wal}. However, this
is clearly not possible when $G$ has an outward pointing cusp, see
\cite{AP}. A convergence result of similar form was proved in
\cite{AP} for Bieberbach polynomials in the Bergman kernel method.

\section{Proofs}

\subsection{Proofs of the results from Section 1}

{\bf Proof of Theorem \ref{thm1.1}.} Using the mapping
$\psi:=\phi^{-1}$, we obtain for any $z\in G$ that
\begin{align*}
|\phi(z)-J_{n,p}(z)| &= \left|\int_{\zeta}^z
\left(\phi'(t)-J_{n,p}'(t)\right) dt\right| \\
&= \left|\int_{0}^{\phi(z)} \left(\phi'(\psi(u)) -
J_{n,p}'(\psi(u))\right)\psi'(u)du\right| \\ &\le
\int_{0}^{\phi(z)} \left|\phi'(\psi(u)) - J_{n,p}'(\psi(u))\right|
|\psi'(u)| |du|,
\end{align*}
where the integration is carried over the segment connecting $0$
and $\phi(z)$ in $D_R$. Since $L$ is rectifiable, the function
under the latter integral belongs to the Hardy class $H^1(D_R)$.
Hence we obtain by the Fej\'er-Riesz inequality (cf. \cite[Theorem
3.13]{Du}) that
\begin{align*}
|\phi(z)-J_{n,p}(z)| &\le \frac{1}{2} \int_{|u|=R}
\left|\phi'(\psi(u)) - J_{n,p}'(\psi(u)\right| |\psi'(u)| |du|
\\ &= \frac{1}{2} \int_L \left|\phi'(t) - Q_{n,p}^p(t)\right| |dt|.
\end{align*}
If $p=1$ then we are done. Applying H\"older's inequality for
$p\ge 2$, we have
\begin{align} \label{4.1}
|\phi(z)-J_{n,p}(z)| &\le \frac{1}{2} \int_L
\left|\left(\phi'(t)\right)^{1/p}-Q_{n,p}(t)\right|\,
\left|\sum_{k=0}^{p-1} \left(\phi'(t)\right)^{k/p}
\left(Q_{n,p}(t)\right)^{p-k-1}\right| |dt| \nonumber \\
&\le \frac{1}{2}\
\left\|\left(\phi'\right)^{1/p}-Q_{n,p}\right\|_p \,
\left\|\sum_{k=0}^{p-1} \left(\phi'\right)^{k/p}
\left(Q_{n,p}\right)^{p-k-1}\right\|_q,
\end{align}
where $q=p/(p-1).$ Observe that
\begin{align*}
\left|\sum_{k=0}^{p-1} \left(\phi'(t)\right)^{k/p}
\left(Q_{n,p}(t)\right)^{p-k-1}\right| &\le \sum_{k=0}^{p-1}
\left|\phi'(t)\right|^{k/p} \left|Q_{n,p}(t)\right|^{p-k-1} \\
&\le \left(\left|\phi'(t)\right|^{1/p}+|Q_{n,p}(t)|\right)^{p-1},
\end{align*}
so that
\begin{align*}
\left\|\sum_{k=0}^{p-1} \left(\phi'\right)^{k/p}
\left(Q_{n,p}\right)^{p-k-1}\right\|_q &\le \left( \int_L
\left(\left|\phi'(t)\right|^{1/p}+|Q_{n,p}(t)|\right)^p |dt|
\right)^{\frac{p-1}{p}} \\ &\le \left(\left\|
\left(\phi'\right)^{1/p}\right\|_p + \|Q_{n,p}\|_p\right)^{p-1},
\end{align*}
by Minkowski's inequality. Since $\left\|
\left(\phi'\right)^{1/p}\right\|_p = \left( \int_L
\left|\phi'(t)\right|\, |dt| \right)^{1/p} = (2\pi R)^{1/p},$ it
follows from \eqref{4.1} that
\begin{align} \label{4.2}
\|\phi-J_{n,p}\|_{\infty} \le \frac{1}{2}\left((2\pi R)^{1/p} +
\|Q_{n,p}\|_p \right)^{p-1}
\left\|(\phi')^{1/p}-Q_{n,p}\right\|_p.
\end{align}
Recall that $\|Q_{n,p}\|_p \le \|1\|_p = l^{1/p}$, by the
definition of $Q_{n,p}$, so that the inequality is proved. The
second statement now follows from Part {\it 2} of Theorem KL.
\\ \db

\noindent{\bf Proof of Theorem \ref{thm1.2}.} Repeating all steps
of the proof of Theorem \ref{thm1.1} up to \eqref{4.2}, but with
$\tilde J_{n,p}$ and $\tilde Q_{n,p}$ instead of $J_{n,p}$ and
$Q_{n,p}$, we obtain the inequality
\[
\|\phi-\tilde J_{n,p}\|_{\infty} \le \frac{1}{2}\left((2\pi
R)^{1/p} + \|\tilde Q_{n,p}\|_p \right)^{p-1}
\left\|(\phi')^{1/p}-\tilde Q_{n,p}\right\|_p.
\]
The proof of the desired inequality is finished by estimating
\begin{align*}
\|\tilde Q_{n,p}\|_p &\le \left\|\tilde Q_{n,p} -
(\phi')^{1/p}\right\|_p + \left\|(\phi')^{1/p}\right\|_p \le
\left\|1-(\phi')^{1/p}\right\|_p + \left\|(\phi')^{1/p}\right\|_p
\\ &\le \|1\|_p + 2 \left\|(\phi')^{1/p}\right\|_p =
l^{1/p} + 2(2\pi R)^{1/p}.
\end{align*}
We also have from the definition of $\tilde Q_{n,p}$ in
\eqref{1.3} that
\[
\left\|(\phi')^{1/p}-\tilde Q_{n,p}\right\|_p \le
\left\|(\phi')^{1/p}-Q_{n,p}\right\|_p,
\]
which tends to 0 as $n\to\infty$ in Smirnov domains, by Part {\it
2} of Theorem KL.\\ \db

We connect the analyticity of $\phi$ on $\overline G$ and the
asymptotics for the leading coefficients of $\tilde Q_{n,p}$ in
the following lemma.
\begin{lemma} \label{lem4.1}
Let $G$ be a Smirnov domain. Set $\tilde Q_{n,p}(z)=\tilde a_{n,p}
z^n + \ldots + \tilde a_{0,p},\ n\in\N.$ If $\phi$ is not analytic
on $\overline G$, then
\begin{align} \label{4.3}
\limsup_{n\to\infty} |\tilde a_{n,p}|^{1/n} =
\frac{1}{\textup{cap}(\overline G)},
\end{align}
where $\textup{cap}(\overline G)$ is the logarithmic capacity of
$\overline G$.
\end{lemma}
\noindent{\bf Proof.} The idea of this proof is suggested by Blatt
and Saff \cite{BS}. We first note that, for any polynomial
$P_n(z)=a_nz^n+\ldots,$ the following holds true
\begin{align} \label{4.4}
|a_{n}| \le \left(\textup{cap}(\overline G)\right)^{-n}
\|P_n\|_{\infty},
\end{align}
see Lemma 4.1 in \cite{BS}. Indeed, if $\Phi$ is the conformal
mapping of $\Omega:=\overline\C\setminus \overline G$ onto the
exterior of the unit disk, normalized by $\Phi(\infty)=\infty$ and
$\lim_{z\to\infty} \Phi(z)/z=1/\textup{cap}(\overline G),$ then
\[
\left|\frac{P_n(z)}{\Phi^n(z)}\right| \le \|P_n\|_{\infty}, \qquad
z\in\Omega,
\]
by the maximum modulus principle for $P_n(z)/\Phi^n(z)$ in
$\Omega.$ Now let $z\to\infty$ to obtain \eqref{4.4}. It follows
from Theorem 1.1 of \cite{Pr97} that
\begin{align} \label{4.5}
\|P_n\|_{\infty} \le c_1 n^{2/p} \|P_n\|_p,
\end{align}
where $c_1>0$ is independent of $n$. Therefore,
\[
\limsup_{n\to\infty} |\tilde a_{n,p}|^{1/n} \le
\frac{\limsup_{n\to\infty} \|\tilde
Q_{n,p}\|_{\infty}^{1/n}}{\textup{cap}(\overline G)} \le
\frac{\limsup_{n\to\infty} \|\tilde
Q_{n,p}\|_p^{1/n}}{\textup{cap}(\overline G)} =
\frac{1}{\textup{cap}(\overline G)},
\]
because of \eqref{4.4}, \eqref{4.5} and $\lim_{n\to\infty}
\|\tilde Q_{n,p}\|_p = \|\left(\phi'\right)^{1/p}\|_p.$ We assume
that
\[
\limsup_{n\to\infty} |\tilde a_{n,p}|^{1/n} <
\frac{1}{\textup{cap}(\overline G)},
\]
and show this leads to a contradiction. Consider a sequence of
Fekete polynomials $F_n,\ n\in\N,$ for $\overline G$, so that
\[
\lim_{n\to\infty} \|F_n\|_{\infty}^{1/n}=\textup{cap}(\overline
G),
\]
see \cite[Sect. 5.5]{Ra}. We define a new sequence $q_n(z):=\tilde
a_{n,p} (z-\zeta) F_{n-1}(z) = \tilde a_{n,p} z^n + \ldots,\
n\in\N.$ It follows from the extremal property \eqref{1.3} that
\begin{align*}
\left\|(\phi')^{1/p}-\tilde Q_{n-1,p}\right\|_p \le
\left\|(\phi')^{1/p} - \left(\tilde Q_{n,p}-q_n\right)\right\|_p
\le \left\|(\phi')^{1/p}-\tilde Q_{n,p}\right\|_p + \|q_n\|_p,
\end{align*}
for $p\in[1,\infty).$ Thus we obtain from the above that
\begin{align*}
\limsup_{n\to\infty} \left(\left\|(\phi')^{1/p}-\tilde
Q_{n-1,p}\right\|_p - \left\|(\phi')^{1/p}-\tilde
Q_{n,p}\right\|_p\right)^{1/n} &\le \limsup_{n\to\infty}
\|q_n\|_p^{1/n} \\ \le \limsup_{n\to\infty} |\tilde a_{n,p}|^{1/n}
\lim_{n\to\infty} \|F_{n-1}\|_{\infty}^{1/n} &< 1.
\end{align*}
Consequently,
\begin{align} \label{4.6}
d:=\limsup_{n\to\infty} \left\|(\phi')^{1/p}-\tilde
Q_{n,p}\right\|_p^{1/n} < 1,
\end{align}
as $\lim_{n\to\infty} \left\|(\phi')^{1/p}-\tilde
Q_{n,p}\right\|_p = 0$ by Theorem KL. If $p\in(0,1)$ then we have
that
\[
\left\|(\phi')^{1/p}-\tilde Q_{n-1,p}\right\|_p^p \le
\left\|(\phi')^{1/p}-\tilde Q_{n,p}\right\|_p^p + \|q_n\|_p^p,
\]
and \eqref{4.6} follows by a similar argument.

Since $\tilde Q_{n,p}$ converges to $(\phi')^{1/p}$ locally
uniformly in $G$ by Theorem KL, we can write
\[
\left(\phi'(z)\right)^{1/p}-\tilde Q_{n,p}(z) =
\sum_{k=1}^{\infty} \left(\tilde Q_{(k+1)n,p}(z)-\tilde
Q_{kn,p}(z)\right), \qquad z\in G.
\]
Thus
\begin{align*}
\left\|(\phi')^{1/p}-\tilde Q_{n,p}\right\|_{\infty} &\le
\sum_{k=1}^{\infty} \left\|\tilde Q_{(k+1)n,p}-\tilde
Q_{kn,p}\right\|_{\infty} \\ &\le \sum_{k=1}^{\infty} c_1
\left((k+1)n\right)^{2/p} \left\|\tilde Q_{(k+1)n,p}-\tilde
Q_{kn,p}\right\|_p,
\end{align*}
by \eqref{4.5}. It is clear from \eqref{4.6} that there exist
$c_2,\varepsilon>0$ such that $d+\varepsilon<1$ and
\[
\left\|(\phi')^{1/p}-\tilde Q_{n,p}\right\|_p < c_2
(d+\varepsilon/2)^n, \qquad n\in\N.
\]
Hence for $p\in[1,\infty)$
\begin{align*}
\left\|\tilde Q_{(k+1)n,p}-\tilde Q_{kn,p}\right\|_p &\le
\left\|(\phi')^{1/p}-\tilde Q_{(k+1)n,p}\right\|_p +
\left\|(\phi')^{1/p}-\tilde Q_{kn,p}\right\|_p \\ &\le 2
\left\|(\phi')^{1/p}-\tilde Q_{kn,p}\right\|_p \le 2 c_2
(d+\varepsilon/2)^{kn},
\end{align*}
and for $p\in(0,1)$
\begin{align*}
\left\|\tilde Q_{(k+1)n,p}-\tilde Q_{kn,p}\right\|_p^p &\le
\left\|(\phi')^{1/p}-\tilde Q_{(k+1)n,p}\right\|_p^p +
\left\|(\phi')^{1/p}-\tilde Q_{kn,p}\right\|_p^p \\ &\le 2
\left\|(\phi')^{1/p}-\tilde Q_{kn,p}\right\|_p^p \le 2 c_2^p
(d+\varepsilon/2)^{knp}.
\end{align*}
It now follows that
\begin{align*}
\left\|(\phi')^{1/p}-\tilde Q_{n,p}\right\|_{\infty} &\le c_3
\sum_{k=1}^{\infty} \left((k+1)n\right)^{2/p}
(d+\varepsilon/2)^{kn} \le c_4 \sum_{k=1}^{\infty}
(d+\varepsilon)^{kn} \\ &\le c_5 (d+\varepsilon)^n,\quad n\in\N.
\end{align*}
The latter estimate is well known to imply that $(\phi')^{1/p}$ is
analytic on $\overline G$, see \cite{Wal}. Contradiction. \\ \db

\noindent{\bf Proof of Theorem \ref{thm1.3}.} Consider the monic
polynomials $\tilde Q_{n,p}(z)/\tilde a_{n,p},\ n\in\N$. Let
$\tilde N$ be a subsequence such that \eqref{4.3} holds along
$\tilde N$ as a regular limit. Then we obtain with the help of
\eqref{4.5} that
\[
\limsup_{\stackrel{n\to\infty}{n\in\tilde N}} \left\|\tilde
Q_{n,p}/\tilde a_{n,p}\right\|_{\infty}^{1/n} =
\lim_{\stackrel{n\to\infty}{n\in\tilde N}} \left|\tilde
a_{n,p}\right|^{-1/n} \limsup_{\stackrel{n\to\infty}{n\in\tilde
N}} \left\|\tilde Q_{n,p}\right\|_{\infty}^{1/n} \le
\textup{cap}(\overline G),
\]
as $\lim_{n\to\infty} \|\tilde Q_{n,p}\|_p =
\|\left(\phi'\right)^{1/p}\|_p.$ Since $\tilde Q_{n,p}$ converge
to $(\phi')^{1/p}$ locally uniformly in $G$, we have
\[
\lim_{n\to\infty} \tilde \nu_{n,p}(E) = 0
\]
for any compact $E\subset G,$ by Hurwitz's theorem. Theorem
\ref{thm1.3} now follows from Theorem 2.1 of \cite{BSS} (or
Theorem 2.1.7 of \cite{AB}). \\
\db

\subsection{Proofs of the results from Section 2}

\noindent{\bf Proof of Theorem \ref{thm2.1}.} We extend the ideas
of \cite{Pr} for $p=2$ to the general case. Let us continue the
mapping $\phi$ conformally beyond the boundary $L$, by using
reflections across the analytic arcs $L_i,\ L=\cup_{i=1}^m L_i.$
Suppose that $\tau_i$ is a mapping such that $L_i=\tau_i([0,1])$,
which is conformal in an open neighborhood of $[0,1].$ Then we can
find a symmetric lens shaped domain $S_i$, bounded by two circular
arcs subtended by $[0,1]$, whose closure is contained in this open
neighborhood of $[0,1].$ Defining
\[
\tilde G := G \cup \left( \cup_{i=1}^m \tau_i(S_i) \right),
\]
we extend $\phi$ into $\tilde G$ as follows:
\[
\phi (z) :=\frac{R^2}{\overline{\phi \left[ \tau_i \left(
\overline{\tau^{-1}_i(z)} \right) \right]} }, \qquad z \in
\tau_i(S_i) \backslash \overline{G},
\]
where $i=1,\ldots,m.$ The boundary $\partial\tilde G$ consists of
$m$ analytic arcs $\Gamma_i$ that share endpoints with the arcs
$L_i$ of $\partial G$:
\[
\partial\tilde G \cap \partial G = \{z_i\}_{i=1}^m,
\]
which are clearly the corner points of $\partial G.$ Since each
$\tau_i,\ i=1,\ldots,m,$ is conformal and has bounded derivative
(together with its inverse) on $S_i$, we obtain the inequalities
\begin{align} \label{4.7}
{\rm dist}(z,\partial G) \ge c_1 \min_{1\le i\le m} |z-z_i|,
\qquad z\in \partial\tilde G,
\end{align}
where ${\rm dist}(z,\partial G)$ is the distance from $z$ to
$\partial G$, and
\begin{align} \label{4.8}
|\gamma| \le c_2 |z-t|, \qquad z,t\in
\partial\tilde G,
\end{align}
where $|\gamma|$ is the length of the shorter arc $\gamma\subset
\partial\tilde G,$ connecting $z$ and $t$. We denote various
positive constants by $c_1,c_2,$ etc.

Let $\Gamma_j$ be an arc of $\partial \tilde G$ with the endpoints
$z_j$ and $z_{j+1}$, and let $\zeta_j \in \Gamma_j$ be a fixed
point, $j=1,\ldots,m.$ Note that $\zeta_j$ divides $\Gamma_j$ into
$\Gamma_j^1$ and $\Gamma_j^2$, so that $\partial \tilde G =
\bigcup_{j =1}^m \bigcup_{i =1}^2 \Gamma_j^i$. We obtain from
Cauchy's integral formula for the continuation of $(\phi')^{1/p}$
into $\tilde G$ that
\begin{align} \label{4.9}
\left(\phi'(z)\right)^{1/p} = \frac{1}{2 \pi i}  \int_{\partial
\tilde G} \frac{\left(\phi'(t)\right)^{1/p}}{t-z}\, dt =
\frac{1}{2 \pi i} \sum_{j =1}^m \sum_{i=1}^2 \int_{\Gamma_j^i}
\frac{\left(\phi'(t)\right)^{1/p}}{t-z}\, dt, \quad z \in \tilde
G.
\end{align}
Hence we need to approximate the functions of the form
\begin{align} \label{4.10}
g(z) := \int_{\gamma} \frac{\left(\phi'(t)\right)^{1/p}}{t-z}\, dt
\end{align}
in $E_p(G)$ norm, where $\gamma$ is any of the arcs $\Gamma_j^i,$
with $i=1,2$ and $j=1,\ldots,m.$

Let $\Omega:=\overline{\C}\setminus\overline{G}.$ Consider the
standard conformal mapping $\Phi:\Omega\to\Delta,$ where
$\Delta:=\{w:|w|>1\},$ normalized by $\Phi(\infty)=\infty$ and
$\Phi'(\infty)>0.$ We define the level curves of $\Phi$ by
\[
L_n:=\{z: |\Phi(z)|=1+1/n\}, \qquad n\in\N.
\]
Denote by $\gamma_2$ the part of $\gamma$ from its endpoint
$\zeta_j\in\Gamma_j$ to the first point $\xi$ of intersection with
$L_n$, so that $\gamma_2 \subset \{z: |\Phi(z)|>1+1/n\}.$ Then
$\gamma_1:=\gamma\setminus\gamma_2$ connects $\xi$ with the corner
point $z_j$ of $L.$ Write
\begin{align} \label{4.11}
g(z) := \int_{\gamma_1}
\frac{\left(\varphi'(t)\right)^{1/p}}{t-z}\, dt + \int_{\gamma_2}
\frac{\left(\varphi'(t)\right)^{1/p}}{t-z}\, dt =: g_1(z) +
g_2(z).
\end{align}
We show that $\|g_1\|_p\to 0$ sufficiently fast as $n\to\infty,$
while $g_2$ is well approximated by polynomials of degree $n.$ To
estimate the norm of $g_1$, we need to know the behavior of
$\left(\phi'\right)^{1/p}$ near the corner point $z_j\in L.$ This
is found from the asymptotic expansion of Lehman \cite{Le}. Assume
that $z_j=0$ and that $\lambda_j\pi,\ 0<\lambda_j<2,$ is the
exterior angle formed by $L$ at this point. Then we have in a
neighborhood of $z_j=0$ that
\[
\phi(z)-\phi(0) = b\ z^{\frac{1}{2-\lambda_j}} +
o\left(z^{\frac{1}{2-\lambda_j}}\right) \qquad \mbox{as } z \to 0,
\]
where $b\neq 0$, and
\[
\phi'(z) = \frac{b}{2-\lambda_j}\ z^{\frac{1}{2-\lambda_j}-1} +
o\left(z^{\frac{1}{2-\lambda_j}-1}\right) \qquad \mbox{as } z \to
0.
\]
Hence there exists a constant $c_3>0$ such that
\begin{align} \label{4.12}
\left|\phi'(z)\right|^{1/p} \le c_3\ |z|^{\alpha}, \qquad z\in
\tilde G \cup \partial \tilde G,
\end{align}
where we set
\[
\alpha:=\frac{1}{p(2-\lambda_j)}-\frac{1}{p}.
\]
For the endpoints $\xi\in L_n$ and $0$ of $\gamma_1$, we let
\[
d_n:=|\xi-0|=|\xi|.
\]
It follows from \eqref{4.8} that
\[
|\gamma_1| \le c_2 d_n.
\]
We now estimate that
\begin{align} \label{4.13}
\|g_1\|_p^p = \int_L \left| \int_{\gamma_1}
\frac{\left(\varphi'(t)\right)^{1/p}}{t-z}\, dt \right|^p\ |dz|
\le c_4 \int_L \left( \int_{\gamma_1} \frac{|t|^{\alpha}
|dt|}{|t-z|} \right)^p\ |dz|,
\end{align}
by \eqref{4.11} and \eqref{4.12}. Note that if $z\in L$ satisfies
$|z|\ge d_n$, then $|t-z|\sim |z|$ by \eqref{4.7}. Consequently,
\begin{align} \label{4.14}
\int_{L\cap\{|z|\ge d_n\}} \left( \int_{\gamma_1}
\frac{|t|^{\alpha} |dt|}{|t-z|} \right)^p\ |dz| &\le c_5
\int_{L\cap\{|z|\ge d_n\}} \left( \frac{d_n^{\alpha+1}}{|z|}
\right)^p\ |dz|  \\ \nonumber &\le c_6 \left\{
\begin{array}{lll}
d_n^{p\alpha+1},\quad &1<p<\infty,\\
d_n^{\alpha+1} |\log d_n|,\quad &p=1,\\
d_n^{p\alpha+p},\quad &\frac{1-\lambda}{2-\lambda}<p<1,
\end{array}
\right.
\end{align}
because $\alpha+1>0$ (which defines the latter range for $p$). On
the other hand, if $z\in L$ satisfies $|z|\le d_n$, then
$|t-z|\sim |t|+|z|$ by \eqref{4.7}, and we obtain by using
\eqref{4.8} that
\begin{align} \label{4.15}
\int_{L\cap\{|z|\le d_n\}} \left( \int_{\gamma_1}
\frac{|t|^{\alpha} |dt|}{|t-z|} \right)^p\ |dz| &\le c_7
\int_0^{c_8 d_n} \left( \int_0^{c_9 d_n} \frac{s^{\alpha}ds}{s+r}
\right)^p\ dr \\ \nonumber &\le c_7 \int_0^{c_8 d_n} \left(
\int_0^r \frac{s^{\alpha}}{r}\,ds + \int_r^{c_9 d_n}
s^{\alpha-1}ds \right)^p\ dr \\ \nonumber &= c_7 \int_0^{c_8 d_n}
\left( \frac{r^{\alpha}}{\alpha+1} +
\frac{(c_9d_n)^{\alpha}-r^{\alpha}}{\alpha} \right)^p\ dr \\
\nonumber &\le c_{10}\, d_n^{p\alpha+1},
\end{align}
for $\alpha\neq 0$. If $\alpha=0$ then we estimate
\begin{align*}
\int_{L\cap\{|z|\le d_n\}} \left( \int_{\gamma_1}
\frac{|dt|}{|t-z|} \right)^p\ |dz| &= \int_{L\cap\{|z|\le d_n\}}
\left( \int_{\gamma_1} \frac{|t|^{1/2}\, |t|^{-1/2}\, |dt|}{|t-z|}
\right)^p\ |dz| \\&\le (c_2\,d_n)^{p/2} \int_{L\cap\{|z|\le d_n\}}
\left( \int_{\gamma_1} \frac{|t|^{-1/2}\, |dt|}{|t-z|} \right)^p\
|dz|
\\&\le c_2^{p/2}\,d_n^{p/2}\, c_{10}\, d_n^{-p/2+1}= c_2^{p/2}\,c_{10}\, d_n,
\end{align*}
as above. Combining \eqref{4.13}-\eqref{4.15}, we have that
\begin{align} \label{4.16}
\|g_1\|_p \le c_{11}\left\{
\begin{array}{lll}
d_n^{\alpha+1/p},\quad &1<p<\infty,\\
d_n^{\alpha+1} |\log d_n|,\quad &p=1,\\
d_n^{\alpha+1},\quad &\frac{1-\lambda}{2-\lambda}<p<1.
\end{array}
\right.
\nonumber \\ \le c_{11}\left\{
\begin{array}{lll}
d_n^{\frac{1}{p(2-\lambda)}},\quad &1<p<\infty,\\
d_n^{\frac{1}{2-\lambda}} |\log d_n|,\quad &p=1,\\
d_n^{\frac{\lambda-1}{p(2-\lambda)}+1},\quad
&\frac{1-\lambda}{2-\lambda}<p<1,
\end{array}
\right.
\end{align}
where $\lambda=\min_{1\le j \le m} \lambda_j.$

The next step is the construction of approximating polynomials
$P_n$ for $g_2$. This is accomplished by using Dzjadyk's kernels
(see, e.g., \cite{ABD}) of the form
\[
K_n(t,z)=\sum_{i=0}^n a_i(t) z^i, \qquad n\in\N,
\]
which approximate the Cauchy kernel. It was proved in Lemma 5 of
\cite{AG} that a sequence of such kernels can be selected, so that
for any fixed $k\in\N$, and for all $t\in\gamma$ with
$|\Phi(t)|\ge 1+1/n,$ we have
\begin{align} \label{4.17}
\left|\frac{1}{t-z} - K_n(t,z)\right| \le c_{12}\,
\frac{d_n^k}{|t-z|^{k+1}}, \qquad z\in L,
\end{align}
for all sufficiently large $n\in\N.$ In particular, \eqref{4.17}
holds for $t\in\gamma_2.$ Define the polynomials
\[
P_n(z):=\int_{\gamma_2} \left(\phi'(t)\right)^{1/p}\,
K_n(t,z)\,dt,
\]
and estimate
\begin{align*}
\|g_2-P_n\|_p^p &= \int_L \left| \int_{\gamma_2} \left(
\frac{1}{t-z}-K_n(t,z) \right) \left(\phi'(t)\right)^{1/p}\, dt
\right|^p\ |dz|
\\ &\le c_{13} d_n^{kp} \int_L \left( \int_{\gamma_2}
\frac{|t|^{\alpha} |dt|}{|t-z|^{k+1}} \right)^p\ |dz|,
\end{align*}
by \eqref{4.17} and \eqref{4.12}. Observe that $|t-z| \sim
|t|+|z|$ for $t\in\gamma_2.$ Therefore, we have for $k>\alpha+1/p$
that
\begin{align*}
\int_L \left( \int_{\gamma_2} \frac{|t|^{\alpha}
|dt|}{|t-z|^{k+1}} \right)^p |dz| &\le c_{14} \int_0^{c_{15}}
\left( \int_{c_{16} d_n}^{c_{17}} \frac{s^{\alpha}ds}{(s+r)^{k+1}} \right)^p dr \\
&\le c_{14} \int_0^{c_{16} d_n} \left( \int_{c_{16} d_n}^{c_{17}}
s^{\alpha-k-1}\,ds \right)^p dr \\ &+ c_{14} \int_{c_{16}
d_n}^{c_{15}} \left(r^{-k-1} \int_{c_{16} d_n}^r
s^{\alpha}\,ds + \int_r^{c_{17}} s^{\alpha-k-1}\,ds \right)^p dr \\
&\le c_{18}\, d_n^{p(\alpha-k)+1} + c_{19}
\int_{c_{16}d_n}^{c_{15}} r^{p(\alpha-k)}\, dr \\
&\le c_{20}\, d_n^{p(\alpha-k)+1}.
\end{align*}
It follows that
\begin{align} \label{4.18}
\|g_2-P_n\|_p \le c_{21}\, d_n^{\alpha+1/p} \le c_{21}\,
d_n^{\frac{1}{p(2-\lambda)}}.
\end{align}
Combining  \eqref{4.16} and \eqref{4.18}, we obtain
\begin{align} \label{4.19}
\|g-P_n\|_p \le \|g_1\|_p + \|g_2-P_n\|_p \le c_{22}\, \left\{
\begin{array}{ll}
d_n^{\frac{1}{p(2-\lambda)}},\quad &1<p<\infty,\\
d_n^{\frac{1}{2-\lambda}} |\log d_n|,\quad &p=1,
\end{array}
\right.
\end{align}
and
\begin{align} \label{4.20}
\|g-P_n\|_p^p \le \|g_1\|_p^p + \|g_2-P_n\|_p^p \le c_{22}\,
d_n^{\frac{\lambda-1}{2-\lambda}+p},\quad
\frac{1-\lambda}{2-\lambda}<p<1.
\end{align}

Recall that $d_n=|\xi|,$ where $\xi\in L_n \cap \gamma_1.$
Applying the results of \cite{Le} to the conformal mapping
$\Psi:=\Phi^{-1},$ we obtain
\[
z=\Psi(\Phi(z))-\Psi(\Phi(0))=a
\left(\Phi(z)-\Phi(0)\right)^{\lambda_j} +
o\left(\left(\Phi(z)-\Phi(0)\right)^{\lambda_j}\right)
\quad\mbox{as }z\to 0,
\]
where $\lambda_j\pi$ is the exterior angle at $z_j=0,$ and $a\neq
0.$ Thus
\[
d_n=|\xi| \le c_{23} \min_{z\in L_n} |z| \le c_{24}\,
n^{-\lambda_j} \le c_{24}\, n^{-\lambda}, \qquad n\in\N,
\]
and
\begin{align} \label{4.21}
\|g-P_n\|_p \le c_{25}\, \left\{
\begin{array}{lll}
n^{-\frac{\lambda}{p(2-\lambda)}},\quad &1<p<\infty,\\
n^{-\frac{\lambda}{2-\lambda}} \log n,\quad &p=1,\\
n^{-\frac{\lambda(\lambda-1)}{p(2-\lambda)}-\lambda},\quad
&\frac{1-\lambda}{2-\lambda}<p<1,
\end{array}
\right.
\end{align}
where $n\ge 2$, by \eqref{4.19}-\eqref{4.20}. Hence there exists a
sequence of polynomials $Q_n$ such that
\begin{align} \label{4.22}
\|(\phi')^{1/p}-Q_n\|_p \le c_{26}\, \left\{
\begin{array}{lll}
n^{-\frac{\lambda}{p(2-\lambda)}},\quad &1<p<\infty,\\
n^{-\frac{\lambda}{2-\lambda}} \log n,\quad &p=1,\\
n^{-\frac{\lambda(\lambda-1)}{p(2-\lambda)}-\lambda},\quad
&\frac{1-\lambda}{2-\lambda}<p<1,
\end{array}
\right.
\end{align}
for $n\ge 2$. Since $\left|\left((\phi')^{1/p}-Q_n\right)\circ
\psi\right|^p |\psi'|$ is subharmonic in $D_R$, we have
\begin{align*}
|1-Q_n(\zeta)|^p &=
\left|\left(\phi'(\zeta)\right)^{1/p}-Q_n(\zeta)\right|^p =
\left|\left(\phi'(\psi(0))\right)^{1/p}-Q_n(\psi(0))\right|^p
|\psi'(0)| \\ &\le \frac{1}{2\pi R}
\left\|\left(\phi'\right)^{1/p}-Q_n\right\|_p^p.
\end{align*}
Thus \eqref{2.1} follows from \eqref{4.22} and the
extremal property \eqref{1.3} of $\tilde Q_{n,p}$, as
\[
\left\|\left(\phi'\right)^{1/p}-\tilde Q_{n,p}\right\|_p \le
\left\|\left(\phi'\right)^{1/p}-(Q_n-Q_n(\zeta)+1)\right\|_p.
\]

The second part of the theorem, stated in \eqref{2.2}, is a direct
consequence of \eqref{2.1} and Theorem \ref{thm1.2}.\\ \db

\noindent{\bf Proof of Theorem \ref{thm2.2}.} We use a combination
of methods employed in the previous proof and in the proof of
Theorem 2.1 of \cite{AP}. Note that the analytic arcs can only
have a polynomial order of contact at the junction points $z_j,\
j=1,\ldots,m,$ as explained in Remark 2.3 of \cite{AP} and its
proof. Thus we have $x^a$-type outward pointing cusps with some
finite $a>1$. Applying analytic continuation via reflection to
$\phi$, we write the Cauchy integral formula \eqref{4.9} for
$(\phi')^{1/p}$, and again reduce the problem to approximation of
the function $g$ in \eqref{4.10}. The only difference from the
proof of Theorem \ref{thm2.1} is that instead of the lens shaped
domain $S_i$ one has to use a symmetric in real axes domain,
bounded by the arcs of $y=\pm Ax^a$ and $y=\pm A(1-x)^a$, where
$A>0$ is sufficiently small (see \cite{AP} for the details). In
this case, we have
\begin{align} \label{4.23}
{\rm dist}(z,\partial G) \ge c_1 \min_{1\le i\le m} |z-z_i|^a,
\qquad z\in \partial\tilde G,
\end{align}
instead of \eqref{4.7}, by Lemma 4.2 of \cite{AP}.

Let $\Phi:\Omega \to \Delta$ be a conformal map of
$\Omega:={\overline \C} \setminus {\overline G}$ onto
$\Delta:=\{w:|w|>1\}$, satisfying the conditions
$\Phi(\infty)=\infty$ and $\Phi'(\infty)>0.$ Define the level
curves of $\Phi$ by
\[
L_u :=\{z \in {\overline \Omega}: |\Phi(z)|=u\}, \quad u > 1.
\]
Let $G_u:=\textup{Int}\, L_u,\ u>1,$ be the domain bounded by
$L_u.$ Denote $\gamma_1:=\gamma \cap \overline{G}_u$ and
$\gamma_2:=\gamma \setminus \gamma_1$,  so that $\gamma_2$ lies
exterior to $L_u$. Hence the function $g_2$ of \eqref{4.11} is
holomorphic in $G_u$, and is well approximable by polynomials.
Namely, we obtain from Theorem 3 of \cite[p. 145]{SL} that there
exists a sequence of polynomials $\{p_n\}_{n=1}^{\infty}$ such
that
\begin{equation}\label{4.24}
\|g_2-p_n\|_{\infty} \leq c_2 \frac{n}{(v-1)^2} \max_{z \in G_v}
|g_2(z)| \ v^{-n}, \quad n \in \N,
\end{equation}
where $c_2$ is an absolute constant and $1<v<u$. On choosing
$u=1+2n^{-s}$ and $v=1+n^{-s}$, with $s \in (0,1)$, we estimate
\begin{eqnarray*}
\max_{z \in G_v} |g_2(z)| &\leq& \int_{\gamma_2}
\frac{|\phi'(t)|^{1/p}} {|t-z|} |dt| \le \frac{c_3}{\dis \min_{z
\in G_v,\, t \in \gamma_2} |t-z|}
\\ &\le& \frac{c_3}{\textup{dist}(L_u,L_v)},
\end{eqnarray*}
where dist$(L_u,L_v)$ is the distance between $L_u$ and $L_v.$
Note that
\[
\textup{dist}(L_u,L_v) \ge c_4(u-v)^2,
\]
by a result of Loewner (see \cite[p. 61]{ABD}), which implies
\[
\textup{dist}(L_u,L_v) \ge c_4 n^{-2s}.
\]
We conclude that
\[
\max_{z \in G_v} |g_2(z)| \le c_5 n^{2s},
\]
and, using (\ref{4.24}), we obtain that
\begin{equation}\label{4.25}
\|g_2-p_n\|_{\infty} \le c_6 n^{1+4s}(1+n^{-s})^{-n} \le c_7
n^{1+4s}e^{-n^{1-s}}, \quad n \in \N.
\end{equation}

For the companion function $g_1$ of \eqref{4.11}, we estimate
\begin{equation}\label{4.26}
\|g_1\|_{\infty} \le \max_{z \in {\overline G}} \int_{\gamma_1}
\frac{|\phi'(t)|^{1/p}|dt|}{|t-z|}  \le \max_{t \in \gamma_1}
\frac{|\phi'(t)|^{1/p}}{\textup{dist}(t,L)},
\end{equation}
since $|\gamma_1|\to 0$ when $n\to\infty.$ We now show that
$\|g_1\|_{\infty}$ is sufficiently small. Indeed, we have by
Corollary 1.4 of \cite{Po} that
\[
|\phi'(t)| \le c_8 \frac{R-|\phi(t)|}{\textup{dist}(t,L)},\quad
t\in G.
\]
Hence
\[
|\phi'(t)| \le c_9 \frac{R-|\phi(t)|}{\textup{dist}(t,L)} \le
c_{10} \frac{|\phi(t)-\phi(z_j)|}{\textup{dist}(t,L)},\quad t\in
\gamma_1,
\]
where $z_j$ is the endpoint of $\gamma_1$ and the cusp point of
$L$. It follows by Lemmas 4.4 and 4.2 of \cite{AP} that
\begin{align*}
|\phi(t) - \phi(z_j)| &\le c_{11}
\dis\exp\left(-\frac{c_{12}}{|t-z_j|^{a-1}}\right) \\ &\le c_{13}
\exp\left(-c_{14}\, [\textup{dist}(t,L)]^b\right), \quad t\in
\gamma_1,
\end{align*}
where $b<0$. Applying these estimates in \eqref{4.26}, we obtain
that
\begin{eqnarray*}
\|g_1\|_{\infty} \le c_{15} \max_{t \in \gamma_1}
\frac{\exp\left(-c_{14}\, [\textup{dist}(t,L)]^b/p\right)}
{[\textup{dist}(t,L)]^{1+1/p}}.
\end{eqnarray*}
Since the function $x^{-1-1/p}\exp(-cx^b)$, where $c>0$ and $b<0$,
is strictly increasing on an interval $(0,x_0)$, we deduce from
the previous inequality that
\begin{equation}\label{4.27}
\|g_1\|_{\infty} \le c_{15} \frac{\exp\left(-c_{14}\,
[\textup{dist}(t_u,L)]^b/p\right)}
{[\textup{dist}(t_u,L)]^{1+1/p}},
\end{equation}
where $t_u \in L_u$ and $u=1+2n^{-s}$ is sufficiently close to 1.
It is known that $\Psi:=\Phi^{-1}$ is H\"{o}lder continuous on
$\overline{\Delta}$ (see Theorem 3 in \cite{NP}), so that
\[
\textup{dist}(t_u,L) \le c_{16} (u-1)^{\beta} \le c_{17}
n^{-s\beta},
\]
for some $\beta>0.$ Hence we obtain from \eqref{4.27} that
\begin{equation}\label{4.28}
\|g_1\|_{\infty} \le c_{18} n^{(1+1/p)s\beta} \exp\left(-c_{19}\,
n^{-s\beta b}/p\right), \quad  n \in \N.
\end{equation}
Combining (\ref{4.25}) and (\ref{4.28}), we have from \eqref{4.11}
that
\[
\|g-p_n\|_{\infty} \le c_{20} \exp\left(-c_{21} n^r\right), \quad
n \in \N,
\]
where $r \in (0,1)$ is any number satisfying $r < \min(1-s,-s\beta
b).$ Furthermore, this immediately implies that there exists a
sequence of polynomials $\{P_n(z)\}_{n=1}^{\infty}$ such that
\begin{equation}\label{4.29}
\|(\phi')^{1/p}-P_n\|_{\infty} \le c_{22} \exp\left(-c_{21}
n^r\right), \quad n \in \N,
\end{equation}
by (\ref{4.9}). That concludes the proof of \eqref{2.3}, since by
the extremal property \eqref{1.3}
\begin{align*}
\left\|\left(\phi'\right)^{1/p}-\tilde Q_{n,p}\right\|_p &\le
\left\|\left(\phi'\right)^{1/p}-(P_n-P_n(\zeta)+1)\right\|_p \\
&\le l^{\frac{1}{p}}
\left\|\left(\phi'\right)^{1/p}-(P_n-P_n(\zeta)+1)\right\|_{\infty}
\le 2 l^{\frac{1}{p}}
\left\|\left(\phi'\right)^{1/p}-P_n\right\|_{\infty}.
\end{align*}
Equation \eqref{2.4} follows from Theorem \ref{thm1.2} and
\eqref{2.3}. \\ \db

\noindent{\it Igor E. Pritsker}

\noindent{\sc Department of Mathematics\\ 401 Mathematical Sciences\\
Oklahoma State University\\ Stillwater, OK 74078-1058, U.S.A.}

\noindent{\bf email: igor@math.okstate.edu}


\begin{thebibliography}{11}

\bibitem{Ah} L. V. Ahlfors, {\it Two numerical methods in conformal mapping},
Experiments in the computation of conformal maps, pp. 45-52.
National Bureau of Standards Applied Mathematics Series, No. 42.
U. S. Government Printing Office, Washington, D.C., 1955.
\bibitem{ABD} V. V. Andrievskii, V. I. Belyi and V. K. Dzjadyk, Conformal
Invariants in Constructive Theory of Functions of a Complex
Variable, World Federation Publishers, Atlanta, 1995.
\bibitem{AB} V. V. Andrievskii and H.-P. Blatt, Discrepancy of
Signed Measures and Polynomial Approximation, Springer-Verlag, New
York, 2002.
\bibitem{AG} V. V. Andrievskii and D. Gaier, {\it Uniform convergence of
Bieberbach polynomials in domains with piecewise quasianalytic
boundary}, Mitt. Math. Sem. Giessen {\bf 211} (1992), 49-60.
\bibitem{AP} V. V. Andrievskii and I. E. Pritsker, {\it Convergence
of Bieberbach polynomials in domains with interior cusps}, J.
d'Analyse Math. {\bf 82} (2000), 315-332.
\bibitem{BS} H.-P. Blatt and E. B. Saff, { \it Behavior of zeros of
polynomials of near best approximation}, J. Approx. Theory {\bf
46} (1986), 323--344.
\bibitem{BSS} H.-P. Blatt, E. B. Saff and M. Simkani, {\it Jentzsch-Szeg\H{o}
type theorems for the zeros of best approximants}, J. London Math.
Soc. {\bf 38} (1988), 307--316.
\bibitem{Du} P. L. Duren, Theory of $H^p$ Spaces, Dover, New
York, 2000.
\bibitem{Ga64} D. Gaier, Konstruktive Methoden der konformen Abbildung,
Springer-Verlag, Berlin, 1964.
\bibitem{Ga92} D. Gaier, {\it On the convergence of the Bieberbach polynomials in
regions with piecewise analytic boundary}, Arch. Math. {\bf 58}
(1992), 289-305.
\bibitem{Ga98} D. Gaier, {\it Polynomial approximation of conformal maps}, Constr.
Approx. {\bf 14} (1998), 27-40.
\bibitem{Ju} G. Julia, Lecons sur la repr\'esentation conforme des aires
simplement connexes, Gauthier-Villars, Paris, 1931.
\bibitem{Ke} M. V. Keldysh, {\it On a class of extremal polynomials}, Dokl.
Akad. Nauk SSSR {\bf 4} (1936), 163-166. (Russian)
\bibitem{KL35} M. V. Keldysh and M. A. Lavrentiev, {\it On the theory of
conformal mappings}, Dokl. Akad. Nauk SSSR {\bf 1} (1935), 85-87.
(Russian)
\bibitem{KL37} M. V. Keldysh and M. A. Lavrentiev, {\it Sur la repr\'esentation
conforme des domaines limit\'es par des courbes rectifiables},
Ann. Sci. \'Ecole Norm. Sup. {\bf 54} (1937), 1-38.
\bibitem{Le} R. S. Lehman, {\it Development of the mapping function at an analytic
corner}, Pacific J. Math. {\bf 7} (1957), 1437-1449.
\bibitem{NP} R. N\"{a}kki and B. Palka, {\it Lipschitz conditions, b-arcwise
connectedness and conformal mappings}, J. d'Analyse Math. {\bf 42}
(1982/83), 38-50.
\bibitem{Po} Ch. Pommerenke, Boundary  Behaviour of Conformal Maps,
Springer-Verlag, Berlin, 1992.
\bibitem{Pr97} I. E. Pritsker, {\it Comparing norms of polynomials in one and
several variables}, J. Math. Anal. Appl. {\bf 216} (1997),
685-695.
\bibitem{Pr} I. E. Pritsker, {\it Approximation of
conformal mapping via the Szeg\H{o} kernel method}, Comp. Methods
and Function Theory {\bf 3} (2003), 79-94.
\bibitem{Ra} T. Ransford, Potential Theory in the Complex Plane,
Cambridge Univ. Press, Cambridge, 1995.
\bibitem{RW} P. C. Rosenbloom and S. E. Warschawski, {\it
Approximation by polynomials}, in ``Lectures on functions of a
complex variable," Ann Arbor, University of Michigan Press, 1955,
pp. 287-302.
\bibitem{SL} V. I. Smirnov and N. A. Lebedev, Functions
of a Complex Variable: Constructive Theory, MIT Press, Cambridge,
1968.
\bibitem{Sz} G. Szeg\H{o}, Orthogonal Polynomials, Amer. Math. Soc., Providence,
1975.
\bibitem{Wal} J. L. Walsh, Interpolation and Approximation by Rational Functions
in the Complex Domain, Colloquium Publications, Vol. 20, Amer.
Math. Soc., Providence, 1969.
\bibitem{Wa} S. E. Warschawski, {\it Recent results in numerical methods of conformal
mapping}, in ``Proceedings of Symposia in Applied Mathematics.
Vol. VI. Mumerical Analysis," McGraw-Hill Book Company, Inc., New
York, 1956, pp. 219-250.

\end{thebibliography}
\end{document}